\documentclass[12pt]{article}
\usepackage{latexsym, amsfonts, amsmath, amssymb, color, graphicx}
\usepackage{graphicx}

\renewcommand{\text}[1]{\ \mbox{#1}\ }

\newcommand{\qed}{\hskip 5mm \rule{2.5mm}{2.5mm}}
\newcommand{\R}{{\mathbb R}}
\newcommand{\C}{{\mathbb C}}

\newcommand{\proof}{{\em Proof:\ }}
\begin{document}
\newtheorem{thm}{Theorem}[section]
\newtheorem{defs}[thm]{Definition}
\newtheorem{lem}[thm]{Lemma}
\newtheorem{note}[thm]{Note}
\newtheorem{cor}[thm]{Corollary}
\newtheorem{prop}[thm]{Proposition}
\renewcommand{\theequation}{\arabic{section}.\arabic{equation}}
\newcommand{\newsection}[1]{\setcounter{equation}{0} \section{#1}}
\renewcommand{\baselinestretch}{1}
%%%%%%%%%%%%% title %%%%%%%%%%%%%%%%%%%%%%%%%%%%%%%
\title{Inverse scattering on the line with a transfer condition 
      }
%%%%%%%%%%%%%%%%%%%%%%%%%%%%%%%%%%%%%%%%%%%%%%
\author{Sonja Currie  \footnote{ Supported by NRF grant no. IFR2011040100017}\, $^{+}$\\
% Richard J. Emmett\, $^{+}$\\
Marlena Nowaczyk \footnote{ Partially supported by Foundation for Polish Science, Programme Homing 2009/9} \,$^{o}$ \\
Bruce A. Watson \footnote{Supported by NRF grant no. IFR2011032400120} \, $^{+}$ \\ \\
\, $^{+}$ School of Mathematics\\
University of the Witwatersrand\\
Private Bag 3, P O WITS 2050, South Africa \\ \\
\, $^{o}$ AGH University of Science and Technology\\
Faculty of Applied Mathematics\\
al. A. Mickiewicza 30, 30-059 Krakow, Poland}
\maketitle
%%%%%%%%%%%%%%%%%%%% abstract %%%%%%%%%%%%%%%%%%%%%%
\abstract{
The inverse scattering problem for Sturm-Liouville operators on the line with a matrix transfer condition at the origin is considered. We show that the transfer matrix can be reconstructed from
the eigenvalues and reflection coefficient. In addition, for potentials with compact essential support,
we show that the potential can be uniquely reconstructed.}
%%%%%%%%%%%%%%%%%%%%%%%%%%%%%%%%%%%%%%%%%%%%%%
%%%%%%%%%%%%%%% introduction %%%%%%%%%%%%%%%%%%%%%%%%%%
\newsection{Introduction\label{sec-intro}}

In this paper we investigate the inverse scattering problem for the differential equation
\begin{equation}\label{diff}
 \ell y:=-\frac{d^2y}{dx^2} + q(x)y = {\zeta}^2 y,\quad \mbox{on } (-\infty,0)\cup (0,\infty),
\end{equation}
in $L^2(-\infty,0)\bigoplus L^2(0,\infty)=L^2(\R)$ with point transfer condition
\begin{equation}\label{tx-condition}
  \left[\begin{array}{c}y(0^+)\\ y'(0^+)\end{array}\right]=M
   \left[\begin{array}{c}y(0^-)\\ y'(0^-)\end{array}\right].
\end{equation}
Here the entries of $M$ are taken to be real, $q \in L^2(\R)$ is assumed to be real valued and obey 
the growth condition
\begin{equation}\label{potential_condition}
\int^\infty_{-\infty} (1+|x|)|q(x)|dx < \infty.
\end{equation}
Note that (\ref{potential_condition}) gives that $q\in L^1(\R)$. As usual we denote
$\displaystyle{f(0^+):=\lim_{t\downarrow 0} f(t)}$ and
$\displaystyle{f(0^-):=\lim_{t\uparrow 0} f(t)}$.
The operator $L$ in $L^2(\R)$ is defined by 
\begin{equation}\label{opL}
Ly=\ell y
\end{equation} on $\R\setminus \{0\}$ for
$y$ in the domain, $D(L)$,  of $L$ where 
\begin{equation}\label{dom-oper}
D(L) = \{y | y, \ell y \in L^2(\R) , y|^{(j)}_{(-\infty,0)},y|^{(j)}_{(0,\infty)}\in AC, j=0,1, y \mbox{ obeys } (\ref{tx-condition})\}. 
\end{equation}
As $q\in L^2(\R)$, $D(L)$ is independent of $q$ which ensures that $D(L)$ is known a priori for the 
inverse problem.

We will only consider point transfer matrices at the origin and henceforth will refer to 
them as transfer matrices. In the physical context the transfer matrix represents a change of 
medium which affects the incident wave as represented by 
components of the matrix. Our transfer matrices will be real 
constant transfer matrices i.e. all components will be constants. 

Hochstadt and Lieberman, in \cite{hoch}, considered the inverse Sturm-Liouville problem 
of the unique determination of the potential on a given interval from 
the spectrum, the boundary conditions and the potential on half of the interval.
These results were generalized to the case of eigenfunctions having a discontinuity at the mid-point 
of the interval in the famous paper by Hald, \cite{hald}, where, in addition,
it was shown that one boundary condition can also be uniquely recovered. 
This in turn was extended by \cite{willis} to 
the case of two interior discontinuites. Similar techniques were then used by \cite{koba} to give a 
uniqueness proof for the inverse Sturm-Liouville problem on a bounded interval with a symmetric 
potential having two interior jump discontinuities. 

Ramm, \cite{ramm}, discusses inverse scattering and spectral one-dimensional problems 
on the half-line in detail.
 Some of the main topics included are, invertibility of the steps in the Gel'fand-Levitan and 
Marchenko inversion procedures, Krein inverse scattering theory and inverse problems.

It should be noted that in \cite{rost}, Hryniv shows that the potential of a Sturm-Liouville operator 
depends analytically and Lipschitz continuously on the spectral data i.e. two spectra or one spectrum 
and the corresponding norming constants. Since he considers $q\in H^{-1}(0,1)$, this means that 
there could be a discontinuity at an interior point of $(0,1)$. Thus, the inverse problem that Hryniv 
considers could be thought of as a discontinuous Sturm-Liouville problem on a finite interval where 
the transfer condition is of a special form which is less general then the transfer condition which we 
are considering in this paper. In \cite{kes1, kes2, kes3} the authors consider the discontinuous 
Sturm-Liouville operator on a finite interval where the boundary conditions may depend on the 
eigenparameter. 
In \cite{kes1} and \cite{kes3} a transfer condition equivalent to taking
$M=\left[\begin{array}{cc}
\alpha & 0\\
0 & \alpha^{-1 }
\end{array}\right]$
in (\ref{tx-condition}) is used, whereas in \cite{kes2} the transfer condition itself is also 
dependent on the eigenparameter. For each of the various types of problems considered in 
\cite{kes1, kes2, kes3} uniqueness theorems for the solution of inverse problems using the 
Titchmarsh-Weyl function and spectral data are proven.

In this paper we solve the following inverse problem. Given the scattering data, using the 
asymptotics developed in \cite{cenw1}, we provide a reconstruction of the transfer matrix $M$ and 
the scattering coefficients. For the case of the potential having compact essential 
support, given the scattering data, one can determine the Titchmarsh-Weyl 
$m$-function for (\ref{diff}) with separated boundary conditions and transfer condition 
(\ref{tx-condition}),
on $[-S,S]$ where ${\rm ess\, supp}(q)\subset [-S,S]$.
 Consequently the potential can be uniquely reconstructed.

In Section 2 the notation and some basic results are presented.  The refection coefficient
is considered in Section 3.  Attention is restricted to the compact essential support in Section 4,
where the main result is presented.  Section 5 is the Appendix in which the details of the 
asymptotics used in this paper are presented in detail.
%%%%%%%%%%%%%%%%%%%%%%%%%%%%%%%%%%%%%%%%%%%%%%%%%%%%%%%%%%%%%%%%%%%%%%%%%%%%%%%%%%%%%%%%%%%%%%%%%%%%%%%%

\newsection{Preliminaries}

The scattering problem considered in this paper can be treated as two classical half-line problems 
interacting via the matrix transfer condition (\ref{tx-condition}) at the origin. 

The operator eigenvalue problem associated with $L$, of (\ref{opL}), can be reformulated as a 
system eigenvalue problem as follows. Let $y_1(t) = y(t)$, $y_2(t) = y(-t)$ and $ {Y}(t) = 
\left(\begin{array}{c}y_1(t) \\ y_2(t)\end{array}\right)$
and consider the differential operator in $L^2(0,\infty) \oplus L^2(0,\infty)$ given by
\begin{equation}\label{diff_eqn_system}
 {T} {Y} := -\frac{d^2  {Y}}{dx^2} + Q {Y} = \zeta^2  {Y},\end{equation}
where $Q(t) = \left(\begin{array}{cc}q(t) & 0 \\ 0 & q(-t) \end{array}\right)$. 
The domain of $T$ is given by
\begin{equation}\label{domT}
D(T) = \{Y\,|\,Y, TY \in (L^2(0, \infty))^2, Y,Y'\in AC, UY(0) = VY'(0)\}
\end{equation}
where
 $U=\left(\begin{array}{cc}1 & -m_{11}\\ 0 & m_{21} \end{array}\right)$ 
and 
$ V= \left(\begin{array}{cc} 0 & -m_{12}\\ 1 & m_{22} \end{array}\right)$. 
Here $m_{ij}$, for $i,j=1,2$, are the entries of the transfer matrix $M$.
As the norm on $L^2(0,\infty) \oplus L^2(0,\infty)$ we take 
$$\| {Y}\|^2 = \int^\infty_0  {Y}^T   {\overline{Y}} dx.$$
It should be noted that 
$Ly = \zeta^2 y, y\in D(L),$ is equivalent to $TY = \zeta^2 {Y}, Y\in D(T)$. 
The transfer matrix scattering problem can now be posed as
\begin{equation}\label{trans_eqn}
TY = \zeta^2 Y, \quad Y \in D(T).
\end{equation}

For $F,G\in D(T)$, define the Lagrange form 
$S(F,G)  :=
 \langle TF,G\rangle - \langle F,TG\rangle,$
for $F,G\in D(T),$
where 
\begin{equation}\label{norm}
\ \langle F , G\rangle = \int_{0}^{\infty} F(x)^T\bar{G}(x)\,dx. 
\end{equation}

It was shown in \cite[Theorem 3.2]{cenw1} that if $\det M \not= 0$ then the operator $T$ is a self-
adjoint operator if and only if $\det M= 1$, and hence, after rescaling, for any $M$ with 
$\det M>0$, see also \cite{wang}. 

\begin{defs}\label{def_Jost}
\cite[p.297]{chadan} The Jost solutions $f_{+,M}(x,\zeta)$ and $f_{-,M}(x,\zeta)$ are the solutions
 of (\ref{diff}) and (\ref{tx-condition}) with
\begin{equation}
\lim_{x \to \infty} e^{-i\zeta x} f_{+,M}(x,\zeta) = 1 = \lim_{x \to -\infty} e^{i\zeta x} f_{-,M}(x,\zeta). 
\end{equation}
\end{defs}

We can now express the Jost solutions  $f_{+,M}(x,\zeta)$ and $f_{-,M}(x,\zeta)$ to 
(\ref{trans_eqn}) in terms of  the
 classical Jost solutions $f_+(x,\zeta)$ and $f_-(x,\zeta)$ (i.e. when $M=I$)
 by
\begin{equation}\label{F+}
f_{+,M}(x,\zeta) := \left\{\begin{array}{c}f_+(x,\zeta), \quad x > 0 \\
h_1(x,\zeta), \quad x < 0 \end{array}\right.,
\end{equation}
\begin{equation}\label{F-}
f_{-,M}(x,\zeta) := \left\{\begin{array}{c}f_-(x,\zeta), \quad x < 0 \\
h_2(x,\zeta), \quad x > 0 \end{array}\right.,
\end{equation}
where $h_1(x,\zeta)$ and $h_2(x,\zeta)$ are solutions of (\ref{diff}) on $(-\infty, 0)$ and $(0, \infty)$ respectively obeying 
\begin{eqnarray*}
\left(\begin{array}{c}h_1(0^-, \zeta) \\ h'_1(0^-,\zeta)\end{array}\right)&=& M^{-1}\left(\begin{array}{c}f_+(0^+, \zeta) \\ f'_+(0^+,\zeta)\end{array}\right),\\
\left(\begin{array}{c}h_2(0^+, \zeta) \\ h'_2(0^+,\zeta)\end{array}\right)&=&M\left(\begin{array}{c}f_-(0^-, \zeta) \\ f'_-(0^-,\zeta)\end{array}\right).
\end{eqnarray*}

For $M=I$ the existence and asymptotic behaviour of the Jost solutions have been well studied, 
see for example \cite{Freiling, March}. In particular
\begin{equation}\label{4_1f}
 f_+(x,\zeta) = e^{i{\zeta}x} + O\left(\frac{C(x)\rho(x)e^{-{\eta}x}}{1 +
|\zeta|}\right), 
\end{equation}
and
\begin{equation}\label{4_3}
f_-(x,\zeta) = e^{-i{\zeta}x} + O\left(\frac{C(-x)\tilde{\rho}(x)e^{{\eta}x}}{1 +
|\zeta|}\right), 
\end{equation}
as $|x|+|\zeta| \rightarrow \infty$, where $\eta = \Im (\zeta)$. Here $C(x)$ is a non-negative, non-increasing function of $x$ and
\begin{equation}
 \rho(x) = \int_x^{\infty} (1+|\tau |)|q(\tau)|\,d \tau,
\qquad \tilde{\rho}(x) = \int_{-\infty}^x (1+|\tau |)|q(\tau)|\,d \tau .
\end{equation}

% As $q$ is real valued it follows that $\overline{f_\pm(x,\pm \overline{\zeta})}$ 
%and $f_{\pm}(x,\pm \zeta )$ are all solutions of (\ref{diff}), $\zeta\in\C$. 
%The independence of $f_+(x,\xi)$ and $\overline{f_+(x,\xi)} = f_+(x,-\xi)$, 
%for $\xi=\zeta \in\R\backslash \{0\}$, was shown in \cite{cenw1}. 
For $\xi \in \R$, see \cite[Sections 2 and 4]{cenw1}, 
the conjugate Jost solutions take the form
\begin{equation}\label{big_F_conj}
 \overline{f}_{+,M}(x,\xi):= 
   \left\{\begin{array}{cc}\overline{f}_+(x,\xi) = f_+(x,-\xi),\quad & x>0 \\
    \overline{h}_1(x,\xi) = h_1(x,-\xi),\quad & x<0 \end{array}\right.
\end{equation}
which obeys the transfer condition at $x=0$.
Being independent (for $\xi\in\R\backslash\{0\}$), the solutions $f_{+,M}(x,\xi)$ and 
$\overline{f}_{+,M}(x,\xi)$ span the solution space of (\ref{diff}), with (\ref{tx-condition}),
so there exist (unique) coefficients $A(\xi)$ and $B(\xi)$ so that
\begin{equation}\label{Big_A_and_B}
f_{-,M}(x,\xi) = A(\xi)\overline{f}_{+,M}(x,\xi) + B(\xi)f_{+,M}(x,\xi).
\end{equation}
Here $A(\xi)$ and $B(\xi)$ are independent of whether $x>0$ or $x<0$, and
they satisfy the equality
$|A(\xi)|^2 - |B(\xi)|^2 = 1$ for $\xi\in\R\backslash\{0\}$.
The reflection coefficient is defined as
\begin{eqnarray*}
  R(\xi) = \frac{B(\xi)}{A(\xi)}, \mbox{ for } \xi \in \R\backslash\{0\}.
\end{eqnarray*}

From \cite{cenw1} we have that 
\begin{equation}\label{A_asymp}
A(\zeta) = \frac{m_{12}}{2i}\zeta + \frac{m_{11} + m_{22}}{2} + \frac{m_{12}}{2}\int_{-\infty}^{\infty} \cos(-\zeta \tau)q(\tau)e^{i\zeta |\tau|}\,d\tau + O\left(\frac{1}{1+|\zeta|}\right),
\end{equation}
 for large $|\zeta|$ and $\Im (\zeta)\geq 0$.  For $|\xi|$ large, 
\begin{equation}\label{B_asymp}
B(\xi) = -\frac{m_{12}}{2i}\xi + \frac{m_{22} - m_{11}}{2} - \frac{m_{12}}{2}\int_{-\infty}^{\infty} \cos(\xi \tau)q(\tau)e^{-i\xi \tau}\,d\tau + O\left(\frac{1}{1+|\xi|}\right).
\end{equation}

%%%%%%%%%%%%%%%%%%%%%%%%%%%%%%%
\newsection{The reflection coefficient}\label{chap_inv_scat_trans_cond}

In this section, given the reflection coefficient and the eigenvalues we will reconstruct the point
 transfer matrix as well as the coefficients $A(\zeta)$ and $B(\xi)$. Moreover, for three special cases 
of the reflection coefficient $R(\xi)$ we will explicitly find the corresponding transfer matrix and 
$A(\zeta)$. 

%In addition, for $q(x)$ with compact essential support, it is shown that the potential 
%can be uniquely reconstructed.

From \cite[p.175]{Hsieh} we have the following representation result.

\begin{lem}\label{log}
Let $f$ be a function analytic in the upper half-plane obeying 
\begin{itemize}
\item $\zeta(f(\zeta)-1)$ is bounded for $\Im (\zeta) \geq 0,$
\item $f(\zeta)$ is continuous for $\zeta \neq 0$ with $\Im (\zeta) \geq 0,$
\item $f(\zeta)\neq 0$ for $\zeta \neq 0$ with $\Im (\zeta) \geq 0,$
\item $\zeta=0$ is a first order pole of $f(\zeta).$
\end{itemize} 
If $F(\zeta) = \log f(\zeta)$, then for $\zeta \in \C$ with $\Im (\zeta)>0$
$$
F(\zeta)=\frac{1}{2\pi i} \int_{-\infty}^{\infty} \frac{2\Re F(\xi)}{\xi-\zeta}\,d\xi= \frac{1}{2\pi i} \int_{-\infty}^{\infty} \frac{\log|f(\xi)|^2}{\xi-\zeta}\,d\xi.
$$
\end{lem}

 \begin{thm}\label{3.2}
For $m_{12}= 0$, given the scattering data, $\{R(\xi), \eta_1, \dots, \eta_N \},$ where $\eta_1,\dots \eta_N$ are the
eigenvalues of (\ref{opL})-(\ref{dom-oper}), the point transfer 
matrix, $M$, is uniquely determined up to $m_{21}$ and a sign condition. In particular 
$$
m_{22}=\pm\sqrt{\frac{1+C_2}{1-C_2}} \quad \makebox{and} \quad m_{11}=\pm\sqrt{\frac{1-C_2}{1+C_2}}
$$
where $C_2 = \lim_{\xi \rightarrow \infty} R(\xi)$ (and this limit exists). Moreover, 
\begin{equation}\label{A2}
A(\zeta ) = \pm\frac{1}{\sqrt{1-C_2^2}}\prod_{j=1}^{N}\frac{\zeta-i\eta_j}{\zeta+i\eta_j}\exp\left[ 
\frac{1}{2\pi i} \int_{-\infty}^{\infty} \frac{\log \left| 1-C_2^2 \right|-\log(1-|R(\xi)|^2)}{\xi-\zeta}\,d\xi
\right].
\end{equation}
 \label{reconstruct_M2}
 \end{thm}
 
\proof 
From (\ref{A_asymp}) and (\ref{B_asymp}), for $m_{12}=0$, we obtain the following asymptotics
\begin{equation}
A(\zeta)=\frac{m_{11} + m_{22}}{2} + O\left(\frac{1}{1+|\zeta|}\right),
\label{A_asymp0}
\end{equation}
 \begin{equation}
B(\xi)= \frac{m_{22} - m_{11}}{2} + O\left(\frac{1}{1+|\xi|}\right).
\label{B_asymp0}
\end{equation}
As ${\rm det} M=1$ and $m_{12}=0$ it follows that 
$m_{11}m_{22}=1$, that is
\begin{eqnarray}
 m_{11}=\frac{1}{m_{22}}, \label{plus-1}
\end{eqnarray} 
making $m_{11}+m_{22}\not= 0$. 
Thus, given $R(\xi)$ and that $m_{12}=0$, we have
\begin{eqnarray*}
R(\xi)= \frac{B(\xi)}{A(\xi)}&=&\frac{\frac{m_{22} - m_{11}}{2} + O\left(\frac{1}{1+|\xi|}\right)}{\frac{m_{11} + m_{22}}{2} + O\left(\frac{1}{1+|\xi|}\right)}\\
&=& \frac{m_{22}-m_{11}}{m_{11}+m_{22}} + O\left(\frac{1}{1+|\xi|}\right)\\
&\to&\frac{m_{22}-m_{11}}{m_{11}+m_{22}}
\end{eqnarray*}
as $|\xi|\to\infty$.
Denote $\displaystyle C_2=\frac{m_{22}-m_{11}}{m_{11}+m_{22}}$.
 then
$$ 
R(\xi)=C_2 + O\left(\frac{1}{1+|\xi|}\right).
$$
Hence $\lim_{\xi \rightarrow \infty} R(\xi)$ exists and is $C_2$.
Since the reflection coefficient $R(\xi)$ is given $C_2$ is known.
By combining (\ref{plus-1}) with the definition of $C_2$  we  have
$$
m_{22}=\pm\sqrt{\frac{1+C_2}{1-C_2}} 
\quad \makebox{and} \quad m_{11}=\pm\sqrt{\frac{1-C_2}{1+C_2}}.
$$
Hence, in the case when $m_{12}=0$,  the point transfer matrix can 
determined up to one parameter ($m_{21}$ being undetermined) from the scattering data. 

We now turn our attention to $A(\zeta)$. Let
\begin{equation}
f(\zeta)=\frac{2}{m_{11}+m_{22}} A(\zeta)\prod_{j=1}^{N}\frac{\zeta+i\eta_j}{\zeta-i\eta_j},
\label{f_prod0}
\end{equation}
where $0<\eta_1<\dots<\eta_N$ and $-\eta_j^2, j=1,\dots,N$ are the eigenvalues 
of (\ref{opL}), (\ref{dom-oper}),  see \cite[Theorem 3.3]{cenw1}.
Using (\ref{A_asymp0}) we get the asymptotic expression $f(\zeta)=1+O(\frac{1}{\zeta}).$ 
Moreover all the properties given in Lemma \ref{log} are obeyed. Thus setting 
$F(\zeta) = \log f(\zeta)$, for $\zeta \in \C$ with $\Im (\zeta)>0$, we have
\begin{eqnarray}\label{star-aug}
F(\zeta)=\frac{1}{2\pi i} \int_{-\infty}^{\infty} \frac{2\Re F(\xi)}{\xi-\zeta}\,d\xi= \frac{1}{2\pi i} \int_{-\infty}^{\infty} \frac{\log|f(\xi)|^2}{\xi-\zeta}\,d\xi.
\end{eqnarray}
Here, as $\xi \in \R$,
\begin{eqnarray}\label{plus-2}
\log|f(\xi)|^2 
&=& 2 \log \left| \frac{2}{m_{11}+m_{22}} \right|+\log|A(\xi)|^2.
\end{eqnarray}
By \cite[Lemma 4.2]{cenw1},  $$|A(\xi)|^2 - |B(\xi)|^2 = 1.$$ 
So, in terms of $R(\xi)$ we have
\begin{equation}\label{ide}
|A(\xi)|^2(1-|R(\xi)|^2)=1,
\end{equation}
which together with (\ref{plus-2}) gives
\begin{eqnarray}
\log|f(\xi)|^2= 2\log \left| \frac{2}{m_{11}+m_{22}} \right|-\log(1-|R(\xi)|^2).
\label{times}
\end{eqnarray}
Combining (\ref{times}) and (\ref{star-aug}) we obtain
$$
F(\zeta)= \frac{1}{2\pi i} \int_{-\infty}^{\infty} \frac{2\log \left| \frac{2}{m_{11}+m_{22}} \right|-\log(1-|R(\xi)|^2)}{\xi-\zeta}\,d\xi.
$$
As $f(\zeta)=e^{F(\zeta)}$ we conclude
\begin{eqnarray*}
A(\zeta)&=&\frac{m_{11}+m_{22}}{2}\prod_{j=1}^{N}\frac{\zeta-i\eta_j}{\zeta+i\eta_j}\exp\left[ 
\frac{1}{2\pi i} \int_{-\infty}^{\infty} \frac{2\log \left| \frac{2}{m_{11}+m_{22}} \right|-\log(1-|R(\xi)|^2)}{\xi-\zeta}\,d\xi
\right]
\\
&=& \pm\frac{1}{\sqrt{1-C_2^2}}\prod_{j=1}^{N}\frac{\zeta-i\eta_j}{\zeta+i\eta_j}\exp\left[ 
\frac{1}{2\pi i} \int_{-\infty}^{\infty} \frac{\log \left| 1-C_2^2 \right|-\log(1-|R(\xi)|^2)}{\xi-\zeta}\,d\xi
\right].  \qed 
\end{eqnarray*}

As a consequence of the above theorem we obtain the following corollary.

\begin{cor}\label{special1}
If $R(\xi) =0$ (i.e. the reflectionless case) then the point transfer matrix has the form $M= \left(\begin{array}{cc}
\pm 1&0\\
m_{21} & \pm 1
\end{array}\right)$ and 
\begin{equation}\label{A_reflectionless}
A(\zeta ) = \pm \prod_{j=1}^{N}\frac{\zeta-i\eta_j}{\zeta+i\eta_j}.
\end{equation}
Furthermore, if $m_{12}=0$ and $|R(\xi)| = \pm C_2$ then
$$
m_{22}=\pm\sqrt{\frac{1+C_2}{1-C_2}}, \quad m_{11}=\pm\sqrt{\frac{1-C_2}{1+C_2}}
$$
and
$$A(\zeta ) = \pm\frac{1}{\sqrt{1-C_2^2}}\prod_{j=1}^{N}\frac{\zeta-i\eta_j}{\zeta+i\eta_j}.$$

\end{cor}

\proof Since $R(\xi) =0$ we have that $B(\xi) =0$. Thus, by (\ref{B_asymp}), $m_{12}=0$  
and Theorem \ref{reconstruct_M2} can be applied with $C_2=0$ to give
 $m_{11} = m_{22} = \pm 1$ and
$$A(\zeta ) = \pm \prod_{j=1}^{N}\frac{\zeta-i\eta_j}{\zeta+i\eta_j}.$$

If $m_{12}=0$ and $|R(\xi)| = \pm C_2$ then Theorem  \ref{reconstruct_M2} can be applied
 to give
$$A(\zeta ) = \pm\frac{1}{\sqrt{1-C_2^2}}\prod_{j=1}^{N}\frac{\zeta-i\eta_j}{\zeta+i\eta_j}. \qed$$

We now prove a similar theorem to Theorem \ref{reconstruct_M2} for the case where $m_{12} \not= 0$. 
 
 \begin{thm}\label{3.4}
For $m_{12}\neq 0$, given the scattering data $\{R(\xi), \eta_1, \dots, \eta_N \},$ the coefficients of the point transfer matrix, $M$, obey
$$m_{12}(C_1m_{11}-m_{21})=1$$ 
where
$ R(\xi)=-1+\frac{2i}{\xi}C_1+ O\left(\frac{1}{\xi^2}\right).$ In this case
\begin{eqnarray*}
A(\zeta)&=&\frac{m_{12}\zeta +i(m_{11}+m_{22})}{2 i}\prod_{j=1}^{N}\frac{\zeta-i\eta_j}{\zeta+i\eta_j}\\
&&\times \exp\left[\frac{1}{2\pi i}\int_{-\infty}^{\infty}\frac{2\log \frac{2}{\,\sqrt{\xi^2 m_{12}^2 +(m_{11} + m_{22})^2 }} - \log (1-|R(\xi)|^2)}{\xi - \zeta}\, d\xi\right]. 
\end{eqnarray*}
\label{reconstruct_M}
 \end{thm}
 
\proof As $m_{12}\ne 0$ from (\ref{A_asymp}), (\ref{B_asymp}) and the definition of $R(\xi)$ we get
\begin{eqnarray*}
R(\xi)&=&
\frac{ -\frac{m_{12}}{2i}\xi + \frac{m_{22} - m_{11}}{2} - \frac{m_{12}}{2}\int_{-\infty}^{\infty} \cos(\xi \tau)q(\tau)e^{-i\xi \tau}\,d\tau + O\left(\frac{1}{1+|\xi|}\right)}{\frac{m_{12}}{2i}\xi + \frac{m_{11} + m_{22}}{2} + \frac{m_{12}}{2}\int_{-\infty}^{\infty} \cos(\xi \tau)q(\tau)e^{i\xi |\tau|}\,d\tau + O\left(\frac{1}{1+|\xi|}\right)} \\
&=& \frac{-1+\frac{i(m_{22}-m_{11})}{m_{12}\xi} -\frac{i}{\xi}\int_{-\infty}^{\infty} \cos(\xi \tau)q(\tau)e^{-i\xi \tau}\,d\tau + O\left(\frac{1}{\xi^2}\right) }{1+\frac{i(m_{11}+m_{22})}{m_{12}\xi} + \frac{i}{\xi}\int_{-\infty}^{\infty} \cos(\xi \tau)q(\tau)e^{i\xi |\tau|}\,d\tau + O\left(\frac{1}{\xi^2}\right)} \\
&=& \left(-1+\frac{i(m_{22}-m_{11})}{m_{12}\xi} -\frac{i}{\xi}\int_{-\infty}^{\infty} \cos(\xi \tau)q(\tau)e^{-i\xi \tau}\,d\tau + O\left(\frac{1}{\xi^2}\right) \right)\\
& & \times \left(1-\frac{i(m_{11}+m_{22})}{m_{12}\xi} - \frac{i}{\xi}\int_{-\infty}^{\infty} \cos(\xi \tau)q(\tau)e^{i\xi |\tau|}\,d\tau + O\left(\frac{1}{\xi^2}\right)\right)\\
 &=&  -1+\frac{i(m_{11}+m_{22})}{m_{12}\xi} + \frac{i(m_{22}-m_{11})}{m_{12}\xi} + \frac{i}{\xi}\int_{-\infty}^{\infty} \cos(\xi \tau)q(\tau)e^{i\xi |\tau|}\,d\tau\\
& & - \frac{i}{\xi}\int_{-\infty}^{\infty} \cos(\xi \tau)q(\tau)e^{-i\xi \tau}\,d\tau  +O\left(\frac{1}{\xi^2}\right)
  \\
&=&  -1+\frac{2i m_{22}}{m_{12}\xi} - \frac{1}{\xi}\int_0^\infty q(\tau) \sin (2\xi \tau) \, d\tau + O\left(\frac{1}{\xi^2}\right).
\end{eqnarray*}
Let 
\begin{equation}\label{c1}
C_1= \frac{m_{22}}{m_{12}}
\end{equation}  then 
\begin{equation}\label{dagger}
R(\xi)=-1+\frac{2i}{\xi}C_1 - \frac{1}{\xi}\int_0^\infty q(\tau) \sin (2\xi \tau) \, d\tau+ O\left(\frac{1}{\xi^2}\right).
\end{equation}
Since the reflection coefficient $R(\xi)$ is given and $\int_0^\infty q(\tau) \sin (2\xi \tau) \, d\tau$ tends to $0$ by the Riemann-Lebesgue Lemma, the constant $C_1$ is given by
$$C_1=\lim_{\xi\to\infty} \frac{\xi (R(\xi)+1)}{2i}.$$
In addition we have that 
\begin{equation}\label{det}
\det M=m_{11}m_{22}-m_{12}m_{21}=1 = m_{12}(C_1m_{11}-m_{21}).
\end{equation} 

Now, let
\begin{equation}
f(\zeta)=\frac{2i}{m_{12}\zeta +i(m_{11}+m_{22})} A(\zeta)\prod_{j=1}^{N}\frac{\zeta+i\eta_j}{\zeta-i\eta_j}.
\label{f_prod}
\end{equation}

Substituting (\ref{A_asymp}) into (\ref{f_prod}) results in
\begin{eqnarray*}
f(\zeta) &=& \left(1 + O\left(\frac{1}{\zeta}\right)\right)\prod_{j=1}^{N}\frac{\zeta+i\eta_j}{\zeta-i\eta_j}\\
&=& 1 + O\left(\frac{1}{\zeta}\right).
\end{eqnarray*}
Therefore, for $|\zeta|$ large, $f(\zeta) -1 =O\left(\frac{1}{\zeta}\right) $ and all the conditions required in Lemma \ref{log} are met by $f(\zeta)$. Thus setting $F(\zeta) = \log f(\zeta)$, for $\zeta \in \C$ with $\Im (\zeta)>0$, we can write 
$$
F(\zeta)=\frac{1}{2\pi i}\int_{-\infty}^{\infty}\frac{2 \log |f(\xi)| }{\xi-\zeta}\, d \xi = \frac{1}{2\pi i}\int_{-\infty}^{\infty}\frac{\log|f(\xi)|^2}{\xi-\zeta}\, d \xi.
$$
Now by (\ref{f_prod}), together with (\ref{ide}),
\begin{eqnarray*}
\log|f(\xi)|^2&=&\log\left| \frac{2 A(\xi)}{m_{12}\xi +i(m_{11}+m_{22})} \right|^2 \\
 &=& 2\log\left|\frac{2}{m_{12}\xi +i(m_{11}+m_{22})} \right|+\log|A(\xi)|^2 \\
 &=& 2\log \frac{2}{\,\sqrt{\xi^2 m_{12}^2 +(m_{11} + m_{22})^2 }} - \log (1-|R(\xi)|^2),
\end{eqnarray*}
giving 
$$
F(\zeta)=\frac{1}{2\pi i}\int_{-\infty}^{\infty}\frac{2\log \frac{2}{\,\sqrt{\xi^2 m_{12}^2 +(m_{11} + m_{22})^2 }} - \log (1-|R(\xi)|^2)}{\xi - \zeta}\, d\xi.
$$
Since $f(\zeta)=e^{F(\zeta)}$ using (\ref{f_prod}) we obtain
\begin{eqnarray*}
A(\zeta)&=&\frac{m_{12}\zeta +i(m_{11}+m_{22})}{2 i}\prod_{j=1}^{N}\frac{\zeta-i\eta_j}{\zeta+i\eta_j}\\
&&\times \exp\left[\frac{1}{2\pi i}\int_{-\infty}^{\infty}\frac{2\log \frac{2}{\,\sqrt{\xi^2 m_{12}^2 +(m_{11}+ m_{22})^2 }} - \log (1-|R(\xi)|^2)}{\xi - \zeta}\, d\xi\right]. \qed
\end{eqnarray*}

The following corollary is a direct consequence of Theorem \ref{reconstruct_M} for the case where the exponential term in $A(\zeta)$ reduces to $1$. 

\begin{cor}\label{special2}
For $m_{12}\not= 0$, if
\begin{equation}\label{Rcond}
\frac{4}{\xi^2 m_{12}^2 +(m_{11} + m_{22})^2} = 1-|R(\xi)|^2,
\end{equation}
then 
\[A(\zeta)=\frac{m_{12}\zeta +i(m_{11}+m_{22})}{2 i}\prod_{j=1}^{N}\frac{\zeta-i\eta_j}{\zeta+i\eta_j}.\]
Moreover, the coefficients of the transfer matrix are determined by the equations
$K_1=~m_{12}^2$, $K_2=(m_{11}+m_{22})^2$, $m_{11}m_{22}-m_{21}m_{12} = 1$ and $m_{22} = C_1m_{12}$.
Here $C_1$ is as in (\ref{dagger}) and $K_1$, $K_2$ are known and obey $K_1>0$ and $K_2\geq 0$. This results in four possibilities for the transfer matrix $M$.
\end{cor}

Clearly, in all of the above results in this section, as $R(\xi)$ is given and we can find $A(\xi)$ from the relevant equations, it is possible to obtain $B(\xi)$ since $B(\xi) = R(\xi)A(\xi)$.
 
%%%%%%%%%%%%%%%%%%%%%%%%%%%%%%%%%%%%%%%%%
\newsection{Compact essential support potentials} 

For the remainder of the paper we will assume that the potential $q(x)$ has compact 
essential support, say ${\rm ess\, supp}(q) \subset [-S,S]$ for some $S>0$.

\begin{lem}\label{fin}
 Let ${\rm ess\, supp}(q) \subset [-S,S]$ for some $S>0$. 
 Given the scattering data $\{R(\xi), \eta_1, \dots, \eta_N \},$ the matrix $W(S,\zeta)$  is uniquely determined.
Here 
\begin{equation}
W(x,\zeta)=\left[ \begin{array}{cc}
w_1(x,\zeta) & w_2(x,\zeta) \\
w_1 '(x,\zeta) & w_2 '(x,\zeta)
 \end{array} \right]\quad\mbox{with}\quad W(-S,\zeta) =
 \left[ \begin{array}{cc}
-1 & 0 \\
0 & 1
 \end{array} \right]=:H.
\label{w2w1bc}
\end{equation} 
and $w_1, w_2$ are solutions of (\ref{diff}).
\end{lem}

\proof
As ${\rm ess\, supp}(q)\subset [-S,S]$, for $x\leq -S$ and $\zeta=\xi\in\R$, we have
\[f_{-,M}(x, \xi) = f_- (x, \xi)= e^{-i\xi x} \quad \makebox{and} \quad \overline{f}_{-,M}(x,\xi) = \overline{f}_-(x, \xi) = e^{i\xi x}.\]
By (\ref{Big_A_and_B}), for $x\geq S$,
\[f_{-,M}(x,\xi) = A(\xi)e^{-i\xi x} + B(\xi)e^{i\xi x},\]
and
\[\overline{f}_{-,M}(x,\xi) = \overline{A}(\xi)e^{i\xi x} + \overline{B}(\xi)e^{-i\xi x}.\]

For $x\leq -S$, as $q$ is essentially zero, 
\begin{eqnarray*}
  w_j(x, \xi) = a_jf_{-,M}(x, \xi)+b_j\overline{f}_{-,M}(x, \xi) = a_je^{-i\xi x} + b_je^{i\xi x}
\end{eqnarray*}
and for $x\geq S$, 
\begin{eqnarray*}
w_j(x, \xi) 
= \hat{a}_jf_{-,M}(x, \xi)+ \hat{b}_j\overline{f}_{-,M}(x, \xi) = \hat{a}_je^{-i\xi x} + \hat{b}_je^{i\xi x}.
\end{eqnarray*}
Thus for $j=1,2,$
$$\left(\begin{array}{c}
\hat{a}_j\\
\hat{b}_j
\end{array}\right) = \left(\begin{array}{cc}
A(\xi) & \overline{B}(\xi)\\
B(\xi) & \overline{A}(\xi)
\end{array}\right)\left(\begin{array}{c}
a_j\\
b_j
\end{array}\right).$$
From the initial value $W(-S,\zeta)$ it follows that 
$$\left(\begin{array}{cc} a_1 & a_2\\ b_1 & b_2 \end{array}\right)=
\left(\begin{array}{cc} -\frac{e^{-iS\xi}}{2} &- \frac{1}{2i\xi} e^{-iS \xi}\\
   -\frac{e^{iS\xi}}{2}&  \frac{1}{2i\xi} e^{iS \xi}\end{array}\right)
$$
so 
\begin{eqnarray*}
\left[w_1(S,\xi)\,\, w_2(S,\xi)\right]&=&\left[e^{-i\xi S}\,\, e^{i\xi S}\right]\left(\begin{array}{cc}
\hat{a}_1 & \hat{a}_2\\
\hat{b}_1 & \hat{b}_2
\end{array}\right)\\
& =& \left[e^{-i\xi S}\,\, e^{i\xi S}\right]\left(\begin{array}{cc}
A(\xi) & \overline{B}(\xi)\\
B(\xi) & \overline{A}(\xi)
\end{array}\right)
\left(\begin{array}{cc} -\frac{e^{-iS\xi}}{2} &- \frac{1}{2i\xi} e^{-iS \xi}\\
   -\frac{e^{iS\xi}}{2}&  \frac{1}{2i\xi} e^{iS \xi}\end{array}\right)
\end{eqnarray*}
By Theorems \ref{3.2} and \ref{3.4}, given $R(\xi)$ and $\eta _1, \dots, \eta _N$, we can reconstruct $A(\xi)$ and hence $B(\xi)$ and thus find $w_2(S,\xi)$ and $w_1(S,\xi)$ as above.
\qed

We now use the approach given in \cite{BBW} together with that found in \cite[p. 28]{Freiling} 
in order to prove the unique determination of the potential $q$ from the scattering data.
Let $v$ be the solution of (\ref{diff}) on $[-S,S]$ obeying the transfer condition 
(\ref{tx-condition}) and satisfying the terminal conditions
$v(S) = 0$  and $v'(S) =1$. 
The entries of $W(x,\lambda )$ are entire functions of $\lambda $ and the determinant is the Wronskian of $w_1$ and $w_2$ and thus is equal to $-1$ for all $x$ and $\lambda $.

The Titchmarsh-Weyl m-function of (\ref{diff}) on $[-S,S]$ for double Dirichlet boundary
 conditions $y(-S) = 0 = y(S)$ and the transfer condition (\ref{tx-condition}) 
is that value of $m$ for which 
\begin{equation}\label{mfn}
\psi :=w_1 + mw_2
\end{equation}
obeys the  terminal condition $\psi (S) =0$.
Now 
\[\psi (-S,\lambda ) = w_1(-S,\lambda ) + m(\lambda )w_2(-S,\lambda ) = -1.\]
Let \[\Delta (\lambda ) := \makebox{Wron}[w_2,v] = w_2v' - vw_2' = -v(-S,\lambda ) = w_2(S,\lambda ).\]
The function $\Delta (\lambda )$ is entire in $\lambda $ and the zeros of $\Delta (\lambda )$ are the eigenvalues of (\ref{diff}) with double Dirichlet boundary conditions and the
transfer condition (\ref{tx-condition}).
In addition, $v$ and $\psi $ are linearly dependent and as $\psi (-S,\lambda ) = -1$ we have that $v(x,\lambda ) = -v(-S,\lambda )\psi (x,\lambda ) $. Hence
\[\psi (x,\lambda ) 
= \frac{v(x,\lambda )}{-v(-S,\lambda )}
= \frac{v(x,\lambda )}{\Delta (\lambda )}.\]
If we also define
\[\Psi (x,\lambda ) = \left[\begin{array}{cc}
\psi (x,\lambda ) & w_2(x,\lambda )\\
\psi '(x,\lambda ) & w_2'(x,\lambda )
\end{array}\right]\]
then from (\ref{mfn}) it follows that
\[\Psi (x,\lambda ) = W(x,\lambda )\left[\begin{array}{cc}
1 & 0\\
m(\lambda )& 1
\end{array}\right],\]
for all $x$, and that $\det \Psi = \det W = -1$. 

\begin{thm}\label{potq}
Given the Titchmarsh-Weyl m-function, $m$, to (\ref{diff}) on $[-S,S]$ with double Dirichlet
 boundary conditions and the transfer condition (\ref{tx-condition}) and $\tilde{m}$,
 the Titchmarsh-Weyl m-function for the same problem but with the potential $q$ replaced by
 $\tilde{q}$. If $m = \tilde{m}$ then $q= \tilde{q}$.
\end{thm}

\proof Let tilde ( $\tilde{ }$ ) of any quantity, in what follows, denote the same quantity as previously defined but for the problem with $q$ replaced by $\tilde{q}$. Since $m = \tilde{m}$ the eigenvalues for the problem with potential $q$ coincide with those for the problem with potential $\tilde{q}$. 
Note that since we have self-adjointness the algebraic multiplicity of an eigenvalue equals the geometric multiplicity and in addition all the eigenvalues are simple. 
Thus $\Delta (\lambda )$ and $\tilde{\Delta }(\lambda )$ have the same zeros, all of which are simple.

From the asymptotics given in the Appendix it can be seen that $\Delta $ is of order $\frac{1}{2}$ and similarly $\tilde{\Delta }$ is of order $\frac{1}{2}$. 
Therefore, as $\Delta $ and $\tilde{\Delta }$ are entire functions of order $\frac{1}{2}$ with the same zeros, we have that
\[\Delta =c\tilde{\Delta }.\]
So \[c = \frac{\Delta }{\tilde{\Delta }}\]
and taking the limit as $\lambda$ tends to $-\infty$ gives that $c=1$. Hence $\Delta = \tilde{\Delta }$.

We now proceed as  in \cite{BBW, Freiling}.
For $\Delta (\lambda ) \not= 0$, set 
\[P(x,\lambda ) = \Psi \tilde{\Psi }^{-1}(x,\lambda ) = - \left[\begin{array}{cc}
\psi (x,\lambda ) & w_2(x,\lambda )\\
\psi '(x,\lambda ) & w_2'(x,\lambda )
\end{array}\right] \left[\begin{array}{cc}
\tilde{w_2}' (x,\lambda ) & -\tilde{w_2}(x,\lambda )\\
-\tilde{\psi }'(x,\lambda ) & \tilde{\psi }(x,\lambda )
\end{array}\right].\]

Since $m = \tilde{m}$ we have that
\[\Psi  = W \left[\begin{array}{cc}
1 & 0\\
m & 1
\end{array}\right] \quad \makebox{and} \quad \tilde{\Psi}  = \tilde{W} \left[\begin{array}{cc}
1 & 0\\
m & 1
\end{array}\right] \] 
therefore $P$ has an analytic extension to the entire function
\[ P = W\tilde{W}^{-1} = - \left[\begin{array}{cc}
w_1 (x,\lambda ) & w_2(x,\lambda )\\
w_1'(x,\lambda ) & w_2'(x,\lambda )
\end{array}\right]\left[\begin{array}{cc}
\tilde{w_2}' (x,\lambda ) & -\tilde{w_2}(x,\lambda )\\
-\tilde{w_1}'(x,\lambda ) & \tilde{w_1}(x,\lambda )
\end{array}\right],\]
and $\det P = 1$.
So 
\begin{eqnarray*}
P_{11} &=& -\psi \tilde{w_2}' + w_2 \tilde{\psi }'\\
&=& w_2\psi ' -\psi w_2' + w_2(\tilde{\psi }' - \psi ') - \psi (\tilde{w_2}' - w_2')\\
&=& 1 + \frac{w_2(\tilde{v}' - v') - v(\tilde{w_2}' - w_2')}{\Delta (\lambda )}.
\end{eqnarray*}

Similarly
\begin{eqnarray*}
P_{12} &=& \psi \tilde{w_2} - w_2 \tilde{\psi }\\
&=& \frac{v\tilde{w_2} -w_2\tilde{v}}{\Delta (\lambda )}.
\end{eqnarray*}

Using Theorem \ref{app1} and Theorem \ref{app2} together with the maximum-modulus principle we obtain that
$P_{11} \equiv 1$ and $P_{12} \equiv 0.$

Hence $\Psi(x,\lambda) = \tilde{\Psi }(x,\lambda)$ and $w_2 (x,\lambda)= \tilde{w_2}(x,\lambda)$ giving that $q = \tilde{q}.$ \qed

\begin{thm}\label{qrecon}
If $q$ has bounded essential support, then from the scattering data $\{R(\xi), \eta_1, \dots, \eta_N \},$
 the potential $q$ of the scattering problem on the line with transfer condition at the origin can
 be reconstructed uniquely.
\end{thm}

\proof
Let ${\rm ess\, supp}(q)\subset [-S,S]$ for some $S>0$, then from Lemma \ref{fin},
 given the scattering data, we can find $w_1(S,\xi)$ and $w_2(S,\xi)$.
But $m=-\frac{w_1(S,\xi)}{w_2(S,\xi)}$, so the Titchmarsh-Weyl $m$-function for
(\ref{diff}) on $[-S,S]$ with double Dirichlet boundary conditions and the transfer condition 
(\ref{tx-condition}) is uniquely determined from the scattering data. 
Now applying Theorem \ref{potq} gives that the potential $q$ is uniquely determined by $m$
on $[-S,S]$.

 To show the uniqueness of $q$ on the whole real line, assume we have two different 
potentials $q$ and $\hat{q}$ with compact essential support. Let $S$ be so large that 
$\makebox{ess\, supp}(q)\cup \makebox{ess\, supp} (\hat{q}) \subset [-S,S]$, 
then, as $q(x)$ is unique on $[-S,S]$ we have $q=\hat{q}$ on $[-S,S]$ and thus on $\R$. \qed

%%%%%%%%%%%%%%%%%%%%%%%%%%%%%%%%%%%%%%%%%%%%%%%%%%%%%%%%%%%%%%%%%%%%% 
\section{Appendix}
\label{appendix}
From \cite{levinson} we have the following asymptotics for $-S\leq x<0$:

\begin{equation}\label{w2}
w_2(x,\lambda ) = \frac{\sin \sqrt{\lambda }(x+S)}{\sqrt{\lambda }} + O\left(\frac{e^{|\Im \sqrt{\lambda }|(x+S)}}{\lambda }\right)
\end{equation}
 \begin{equation}\label{w2prime}
w_2'(x,\lambda ) = \cos \sqrt{\lambda }(x+S) + O\left(\frac{e^{|\Im \sqrt{\lambda }|(x+S)}}{\sqrt{\lambda }}\right)
\end{equation}
and for $S\geq x>0$
\begin{equation}\label{v}
v(x,\lambda ) = \frac{-\sin \sqrt{\lambda }(S-x)}{\sqrt{\lambda }} + O\left(\frac{e^{|\Im \sqrt{\lambda }|(S-x)}}{\lambda }\right)
\end{equation}
\begin{equation}\label{vprime}
v'(x,\lambda ) = \cos \sqrt{\lambda }(S-x) + O\left(\frac{e^{|\Im \sqrt{\lambda }|(S-x)}}{\sqrt{\lambda }}\right).
\end{equation}
Since $w_2$ and $v$ obey the transfer condition (\ref{tx-condition}) we obtain the following: 

For $m_{12} \not= 0$ and $m_{22} \not= 0$
\begin{equation}\label{w2xgz}
w_2(0^+,\lambda ) = m_{12} \cos \sqrt{\lambda }S + O\left(\frac{e^{|\Im \sqrt{\lambda }|S}}{\sqrt{\lambda }}\right)
\end{equation}
\begin{equation}\label{w2primexgz}
w_2'(0^+,\lambda ) = m_{22} \cos \sqrt{\lambda }S + O\left(\frac{e^{|\Im \sqrt{\lambda }|S}}{\sqrt{\lambda }}\right).
\end{equation}  

For $m_{12} = 0$ and $m_{22} \not= 0$
\begin{equation}\label{w2xgz2}
w_2(0^+,\lambda ) = m_{11} \frac{\sin \sqrt{\lambda }S}{\sqrt{\lambda }} + O\left(\frac{e^{|\Im \sqrt{\lambda }|S}}{\lambda }\right)
\end{equation}
\begin{equation}\label{w2primexgz2}
w_2'(0^+,\lambda ) = m_{22} \cos \sqrt{\lambda }S + O\left(\frac{e^{|\Im \sqrt{\lambda }|S}}{\sqrt{\lambda }}\right).
\end{equation}  

For $m_{12} \not= 0$ and $m_{22} = 0$
\begin{equation}\label{w2xgz3}
w_2(0^+,\lambda ) = m_{12} \cos \sqrt{\lambda }S + O\left(\frac{e^{|\Im \sqrt{\lambda }|S}}{\sqrt{\lambda }}\right)
\end{equation}
\begin{equation}\label{w2primexgz3}
w_2'(0^+,\lambda ) = m_{21} \frac{\sin \sqrt{\lambda }S}{\sqrt{\lambda }} + O\left(\frac{e^{|\Im \sqrt{\lambda }|S}}{\lambda }\right).
\end{equation}  

For $m_{12}\not= 0 $ and $m_{11}\not= 0$
\begin{equation}\label{vxgz}
v(0^-,\lambda ) = -m_{12} \cos \sqrt{\lambda }S + O\left(\frac{e^{|\Im \sqrt{\lambda }|S}}{\sqrt{\lambda }}\right)
\end{equation}
\begin{equation}\label{vprimexgz}
v'(0^-,\lambda ) = m_{11} \cos \sqrt{\lambda }S + O\left(\frac{e^{|\Im \sqrt{\lambda }|S}}{\sqrt{\lambda }}\right).
\end{equation}  

For $m_{12} = 0 $ and $m_{11}\not= 0$
\begin{equation}\label{vxgz2}
v(0^-,\lambda ) = -m_{22} \frac{\sin \sqrt{\lambda }S}{\sqrt{\lambda }} + O\left(\frac{e^{|\Im \sqrt{\lambda }|S}}{\lambda }\right)
\end{equation}
\begin{equation}\label{vprimexgz2}
v'(0^-,\lambda ) = m_{11} \cos \sqrt{\lambda }S + O\left(\frac{e^{|\Im \sqrt{\lambda }|S}}{\sqrt{\lambda }}\right).
\end{equation}  

For $m_{12}\not= 0 $ and $m_{11}= 0$
\begin{equation}\label{vxgz3}
v(0^-,\lambda ) = -m_{12} \cos \sqrt{\lambda }S + O\left(\frac{e^{|\Im \sqrt{\lambda }|S}}{\sqrt{\lambda }}\right)
\end{equation}
\begin{equation}\label{vprimexgz3}
v'(0^-,\lambda ) = m_{21} \frac{\sin \sqrt{\lambda }S}{\sqrt{\lambda }} + O\left(\frac{e^{|\Im \sqrt{\lambda }|S}}{\lambda }\right).
\end{equation}  

Thus extending $w_2$ we obtain for $S\geq x>0$ that if $m_{12} \not= 0$ then
\begin{equation}\label{w2xgzasy}
w_2(x,\lambda ) = m_{12} \cos \sqrt{\lambda }S\cos \sqrt{\lambda }x + O\left(\frac{e^{|\Im \sqrt{\lambda }(S+x)|}}{\sqrt{\lambda }}\right)
\end{equation}
\begin{equation}\label{w2primexgzasy}
w_2'(x,\lambda ) =- m_{12} (\cos \sqrt{\lambda }S)(\sqrt{\lambda }\sin \sqrt{\lambda }x) + O\left(e^{|\Im \sqrt{\lambda }(S+x)|}\right).
\end{equation}  
If $m_{12} = 0$
\begin{equation}\label{w2xgz2asy}
w_2(x,\lambda ) = m_{11} \frac{\sin \sqrt{\lambda }S}{\sqrt{\lambda }}\cos \sqrt{\lambda }x +  m_{22} \frac{\sin \sqrt{\lambda }x}{\sqrt{\lambda }}\cos \sqrt{\lambda }S  + O\left(\frac{e^{|\Im \sqrt{\lambda }(S+x)|}}{\lambda }\right)
\end{equation}
\begin{equation}\label{w2primexgz2asy}
w_2'(x,\lambda ) =-m_{11}\sin \sqrt{\lambda }S  \sin \sqrt{\lambda }x + m_{22} \cos \sqrt{\lambda }S \cos \sqrt{\lambda }x + O\left(\frac{e^{|\Im \sqrt{\lambda }(S+x)|}}{\sqrt{\lambda }}\right).
\end{equation}  
Similarly we can extend $v$ to obtain for $-S\leq  x<0$ that if $m_{12} \not= 0$ then
\begin{equation}\label{vxgzasy}
v(x,\lambda ) = -m_{12} \cos \sqrt{\lambda }S \cos \sqrt{\lambda }x + O\left(\frac{e^{|\Im \sqrt{\lambda }(S-x)|}}{\sqrt{\lambda }}\right)
\end{equation}
\begin{equation}\label{vprimexgzasy}
v'(x,\lambda ) = m_{12} \sqrt{\lambda } \cos \sqrt{\lambda }S \sin \sqrt{\lambda }x + O\left(e^{|\Im \sqrt{\lambda }(S-x)|}\right).
\end{equation}  
If $m_{12} =0$ then
\begin{equation}\label{vxgz2asy}
v(x,\lambda ) = -m_{22} \frac{\sin \sqrt{\lambda }S}{\sqrt{\lambda }}\cos \sqrt{\lambda }x + m_{11}\frac{\sin \sqrt{\lambda }x}{\sqrt{\lambda }}\cos \sqrt{\lambda }S + O\left(\frac{e^{|\Im \sqrt{\lambda }(S-x)|}}{\lambda }\right)
\end{equation}
\begin{equation}\label{vprimexgz2asy}
v'(x,\lambda ) = m_{22}\sin \sqrt{\lambda }S \sin \sqrt{\lambda }x + m_{11} \cos \sqrt{\lambda }S \cos \sqrt{\lambda }x + O\left(\frac{e^{|\Im \sqrt{\lambda }(S-x)|}}{\sqrt{\lambda }}\right).
\end{equation}  

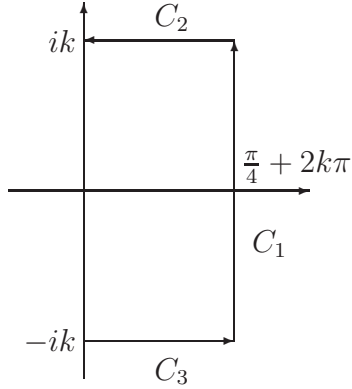
\begin{figure}[h]
\begin{center}
\setlength{\unitlength}{1cm}
\begin{picture}(5,5)
 \put(0,2.5){\vector(1,0){4}}
 \put(1,0){\vector(0,1){5}}
 %\put(1,2.5){\circle{0.3}}
 \put(1,0.5){\vector(1,0){2}}
 \put(3,0.5){\vector(0,1){4}}
 \put(3,4.5){\vector(-1,0){2}}
 \put(3.1,2.6){\makebox(0,0)[bl]{$\frac{\pi}{4}+2k\pi$}}
 \put(0.9,0.5){\makebox(0,0)[r]{$-ik$}}
 \put(0.9,4.5){\makebox(0,0)[r]{$ik$}}
 \put(3.7,1.8){\makebox(0,0)[r]{$C_1$}}
 \put(2.4,4.8){\makebox(0,0)[r]{$C_2$}}
 \put(2.4,0.1){\makebox(0,0)[r]{$C_3$}}
\end{picture}
\end{center}
\caption{$\Gamma_k$ in the $S\sqrt{\lambda}$-plane.}
\label{fig1}
\end{figure}

\begin{lem}\label{deltam120}
For $m_{12}=0$ on the contour $\Gamma_k$ given in Figure \ref{fig1}
\begin{equation}
\frac{1}{|\Delta(\lambda)|}=O(\sqrt{\lambda}e^{-2S|\Im \sqrt{\lambda}|}).
\end{equation}
\end{lem}

\proof 
Note that
$$
\Delta(\lambda)=w_2(S)=\left(m_{11}+\frac{1}{m_{11}}\right)\frac{\sin 2 \sqrt{\lambda}S}{2\sqrt{\lambda}}+O\left(\frac{e^{2S|\Im \sqrt{\lambda}|}}{\lambda} \right).
$$

Consider the contour indicated below for $k\in \mathbb N.$

Let $\displaystyle \lambda=\frac{z^2}{S^2}$ then on $C_1$ 
the variable $
z=i\rho + \frac{\pi}{4}+2k\pi, -k\leq\rho\leq k, $
and
\begin{eqnarray*}
|w_2(S)|&=&\left|\left(m_{11}+\frac{1}{m_{11}}\right)\frac{S}{2}\frac{\sin(2i\rho+\frac{\pi}{2})}{i\rho+\frac{\pi}{4}+2k\pi} \right| +O\left(\frac{e^{2|\rho|}}{k^2} \right)\\
&\geq & \left| \left(m_{11}+\frac{1}{m_{11}}\right)\right| \frac{S}{6k\pi}\,|\cos(2i\rho)|+ O\left(\frac{e^{2|\rho|}}{k^2} \right) \\
& \geq & \left| \left(m_{11}+\frac{1}{m_{11}}\right)\right| \frac{S}{12k\pi}\, e^{2|\rho|} + O\left(\frac{e^{2|\rho|}}{k^2} \right).
\end{eqnarray*}

Therefore 
$$
\frac{1}{|w_2(S)|}\leq \frac{12k\pi}{S}\left| \frac{m_{11}}{1+m_{11}^2} \right| e^{-2|\rho|}+O(e^{-2|\rho|})=O(\sqrt{\lambda}e^{-2S|\Im \sqrt{\lambda}|}).
$$

Considering $C_2$ and $C_3$ let $z=\pm ik+t,$ where $t\in [0,\frac{\pi}{4}+2k\pi]$ then for $|k|$ large
\begin{eqnarray*}
|w_2(S)|&=&\left|\left(m_{11}+\frac{1}{m_{11}}\right)\frac{S}{2}\frac{\sin(\pm 2ik +2t)}{\pm ik+t} \right| + O\left(\frac{e^{2|k|}}{k^2} \right)\\
&=&\left|\left(m_{11}+\frac{1}{m_{11}}\right)\frac{e^{\mp 2k+2it}-e^{\pm 2k -2it}}{2i(\pm ik + t)} \right|\frac{S}{2} + O\left(\frac{e^{2|k|}}{k^2} \right)\\
&\geq & \left| \left(m_{11}+\frac{1}{m_{11}}\right)\right| \frac{S}{2(\pm ik+t)} \frac{e^{2|k|}}{4}+ O\left(\frac{e^{2|k|}}{k^2} \right).\\
\end{eqnarray*}
Thus 
$$
\frac{1}{|w_2(S)|}=O(\sqrt{\lambda} e^{-2S|\Im \sqrt{\lambda}|}). \qed
$$

\begin{figure}[h]
\begin{center}
\setlength{\unitlength}{1cm}
\begin{picture}(5,5)
 \put(0,2.5){\vector(1,0){4}}
 \put(1,0){\vector(0,1){5}}
 %\put(1,2.5){\circle{0.3}}
 \put(1,0.5){\vector(1,0){2}}
 \put(3,0.5){\vector(0,1){4}}
 \put(3,4.5){\vector(-1,0){2}}
 \put(3.1,2.6){\makebox(0,0)[bl]{$2k\pi$}}
 \put(0.9,0.5){\makebox(0,0)[r]{$-ik$}}
 \put(0.9,4.5){\makebox(0,0)[r]{$ik$}}
 \put(3.7,1.8){\makebox(0,0)[r]{$U_1$}}
 \put(2.4,4.8){\makebox(0,0)[r]{$U_2$}}
 \put(2.4,0.1){\makebox(0,0)[r]{$U_3$}}
\end{picture}
\end{center}
\caption{$\Upsilon_k$ in the $S\sqrt{\lambda}$-plane.}
\label{fig2}
\end{figure}
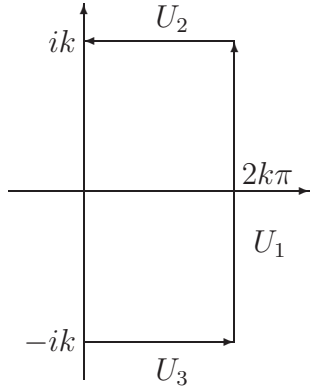

\begin{lem}\label{deltam12n0}
For $m_{12}\neq 0,$ on the contour $\Upsilon_k$ given in Figure \ref{fig2},
\begin{equation}
\frac{1}{|\Delta(\lambda)|}=O(e^{-2S|\Im \sqrt{\lambda}|}).
\end{equation}
\end{lem}

\proof 
In this case
$$
\Delta(\lambda)=w_2(S)=m_{12}\cos^2\sqrt{\lambda}S+O\left(\frac{e^{2S|\Im \sqrt{\lambda}|}}{\sqrt{\lambda}} \right).
$$
Hence we consider the contour indicated below for $k\in \mathbb N.$

Again, let $\displaystyle \lambda=\frac{z^2}{S^2}$ then on $U_1$ the variable
$
z=i\rho + 2k\pi,-k\leq\rho\leq k,
$
and
\begin{eqnarray*}
|w_2(S)|&=& m_{12}\cos ^2 (i\rho +2k\pi) +O\left(\frac{e^{2|\rho|}}{k} \right)\\
&=& m_{12}\left(\frac{e^{-\rho}+e^{\rho}}{2} \right)^2 +O\left(\frac{e^{2|\rho|}}{k} \right)\\
&\geq & m_{12}\frac{e^{2|\rho|}}{4}+ O\left(\frac{e^{2|\rho|}}{k} \right).
\end{eqnarray*}
Therefore 
$$
\frac{1}{|w_2(S)|}= O(e^{-2S|\Im \sqrt{\lambda}|}).
$$

Similarly on $U_2$ and $U_3$ set $z=\pm ik+t$ for $t\in [0, 2k\pi].$ Then for large $|k|$ we have that
\begin{eqnarray*}
|w_2(S)|&=& m_{12}\cos ^2 (\pm ik+t) +O\left(\frac{e^{2|k|}}{k} \right)\\
&=& m_{12}\left(\frac{e^{\mp k+it}+e^{\pm k-it}}{2} \right)^2 +O\left(\frac{e^{2|k|}}{k} \right)\\
&\geq & m_{12}\frac{e^{2|k|}}{4}+ O\left(\frac{e^{2|k|}}{k} \right)
\end{eqnarray*}
giving that
$$
\frac{1}{|w_2(S)|}= O(e^{-2S|\Im \sqrt{\lambda}|}). \qed
$$

\begin{thm}
\label{app1}
On the contours $\Gamma_k$ (for $m_{12}=0$) and $\Upsilon_k$ (for $m_{12}\ne 0$) given in Figure \ref{fig1} and Figure \ref{fig2} respectively
$$
\frac{v\tilde{w_2}-\tilde{v}w_2}{|\Delta(\lambda)|}=O\left(\frac{1}{\lambda}\right),
\quad\mbox{as}\quad k\to\infty.
$$
\end{thm}

\proof
We need to consider four cases. Firstly if $m_{12}=0$ and $x>0$ then
\begin{eqnarray*}
&& v\tilde{w_2}-\tilde{v}w_2 \\
&=& \left( \frac{-\sin \sqrt{\lambda }(S-x)}{\sqrt{\lambda }} + O\left(\frac{e^{|\Im \sqrt{\lambda }|(S-x)}}{\lambda }\right) \right)  \\
& &\times  \left( m_{11} \frac{\sin \sqrt{\lambda }S}{\sqrt{\lambda }}\cos \sqrt{\lambda }x +  m_{22} \frac{\sin \sqrt{\lambda }x}{\sqrt{\lambda }}\cos \sqrt{\lambda }S  + O\left(\frac{e^{|\Im \sqrt{\lambda }(S+x)|}}{\lambda }\right) \right) \\
&-& \left( \frac{-\sin \sqrt{\lambda }(S-x)}{\sqrt{\lambda }} + O\left(\frac{e^{|\Im \sqrt{\lambda }|(S-x)}}{\lambda }\right) \right)  \\
& &\times  \left( m_{11} \frac{\sin \sqrt{\lambda }S}{\sqrt{\lambda }}\cos \sqrt{\lambda }x +  m_{22} \frac{\sin \sqrt{\lambda }x}{\sqrt{\lambda }}\cos \sqrt{\lambda }S  + O\left(\frac{e^{|\Im \sqrt{\lambda }(S+x)|}}{\lambda }\right) \right)\\
&=& O\left(\frac{e^{2S|\Im \sqrt{\lambda}|}}{\lambda ^{3/2}} \right).
\end{eqnarray*}
Hence using Lemma \ref{deltam120} we obtain the result for this particular case.

Next consider $m_{12}=0$ and $x<0$ then
\begin{eqnarray*}
&& v\tilde{w_2}-\tilde{v}w_2 \\
&=& \left(  -m_{22} \frac{\sin \sqrt{\lambda }S}{\sqrt{\lambda }}\cos \sqrt{\lambda }x + m_{11}\frac{\sin \sqrt{\lambda }x}{\sqrt{\lambda }}\cos \sqrt{\lambda }S + O\left(\frac{e^{|\Im \sqrt{\lambda }(S-x)|}}{\lambda }\right) \right) \\
& & \times \left(\frac{\sin \sqrt{\lambda }(x+S)}{\sqrt{\lambda }} + O\left(\frac{e^{|\Im \sqrt{\lambda }|(x+S)}}{\lambda }\right) \right) \\
&-&  \left(  -m_{22} \frac{\sin \sqrt{\lambda }S}{\sqrt{\lambda }}\cos \sqrt{\lambda }x + m_{11}\frac{\sin \sqrt{\lambda }x}{\sqrt{\lambda }}\cos \sqrt{\lambda }S + O\left(\frac{e^{|\Im \sqrt{\lambda }(S-x)|}}{\lambda }\right) \right) \\
& & \times \left(\frac{\sin \sqrt{\lambda }(x+S)}{\sqrt{\lambda }} + O\left(\frac{e^{|\Im \sqrt{\lambda }|(x+S)}}{\lambda }\right) \right) \\
&=& O\left(\frac{e^{2S|\Im \sqrt{\lambda}|}}{\lambda ^{3/2}} \right).
\end{eqnarray*}
Again using Lemma \ref{deltam120} gives the required result for $x<0.$

So far we have shown that for $m_{12}=0$ 
$$ 
\frac{v\tilde{w_2}-\tilde{v}w_2}{|\Delta(\lambda)|}=O\left(\frac{1}{\lambda}\right).
$$

It remains to show that the result holds for $m_{12}\neq 0.$ Again we will consider $x>0$ and $x<0$ separately. Let $x>0$ then
\begin{eqnarray*}
&& v\tilde{w_2}-\tilde{v}w_2 \\
& = & \left( \frac{-\sin \sqrt{\lambda }(S-x)}{\sqrt{\lambda }} + O\left(\frac{e^{|\Im \sqrt{\lambda }|(S-x)}}{\lambda }\right) \right)  \left(  m_{12} \cos \sqrt{\lambda }S\cos \sqrt{\lambda }x + O\left(\frac{e^{|\Im \sqrt{\lambda }(S+x)|}}{\sqrt{\lambda }}\right) \right)\\
&-& \left( \frac{-\sin \sqrt{\lambda }(S-x)}{\sqrt{\lambda }} + O\left(\frac{e^{|\Im \sqrt{\lambda }|(S-x)}}{\lambda }\right) \right)  \left(  m_{12} \cos \sqrt{\lambda }S\cos \sqrt{\lambda }x + O\left(\frac{e^{|\Im \sqrt{\lambda }(S+x)|}}{\sqrt{\lambda }}\right) \right) \\
&=& O\left( \frac{e^{2S|\Im \sqrt{\lambda}|}}{\lambda } \right).
\end{eqnarray*}
The result now follows by Lemma \ref{deltam12n0}.

Lastly if $m_{12}\neq 0$ and $x<0$ then
\begin{eqnarray*}
&& v\tilde{w_2}-\tilde{v}w_2 \\
&=& \left(-m_{12} \cos \sqrt{\lambda }S \cos \sqrt{\lambda }x + O\left(\frac{e^{|\Im \sqrt{\lambda }(S-x)|}}{\sqrt{\lambda }}\right) \right) 
\left(\frac{\sin \sqrt{\lambda }(x+S)}{\sqrt{\lambda }} + O\left(\frac{e^{|\Im \sqrt{\lambda }|(x+S)}}{\lambda }\right) \right) \\
&-& \left(-m_{12} \cos \sqrt{\lambda }S \cos \sqrt{\lambda }x + O\left(\frac{e^{|\Im \sqrt{\lambda }(S-x)|}}{\sqrt{\lambda }}\right) \right) 
\left(\frac{\sin \sqrt{\lambda }(x+S)}{\sqrt{\lambda }} + O\left(\frac{e^{|\Im \sqrt{\lambda }|(x+S)}}{\lambda }\right) \right) \\
&=& O\left( \frac{e^{2S|\Im \sqrt{\lambda}|}}{\lambda } \right).
\end{eqnarray*}
Lemma \ref{deltam12n0} again leads to the required result. This completes all four cases. \qed

Similarly we now prove the following theorem.

\begin{thm}
\label{app2}
On the contours $\Gamma_k$ (for $m_{12}=0$) and $\Upsilon_k$ (for $m_{12}\ne 0$) given in Figure \ref{fig1} and Figure \ref{fig2} respectively
$$
\frac{w_2(\tilde{v}' -v')-v(\tilde{w_2}'-w_2')}{|\Delta(\lambda)|}=O\left(\frac{1}{\sqrt{\lambda}}\right)\quad\mbox{as}\quad k\to\infty.
$$
\end{thm}

\proof
We again need to consider the same four cases as in Theorem \ref{app1}. So we start with $m_{12} = 0$ and $x>0$. In this case
\begin{eqnarray*}
&& w_2(\tilde{v}' -v')-v(\tilde{w_2}'-w_2')\\
&=& \left( m_{11} \frac{\sin \sqrt{\lambda }S}{\sqrt{\lambda }}\cos \sqrt{\lambda }x +  m_{22} \frac{\sin \sqrt{\lambda }x}{\sqrt{\lambda }}\cos \sqrt{\lambda }S  + O\left(\frac{e^{|\Im \sqrt{\lambda }(S+x)|}}{\lambda }\right) \right) \left( O\left(\frac{e^{|\Im \sqrt{\lambda }|(S-x)}}{\sqrt{\lambda }}\right)\right)\\
&-& \left(\frac{-\sin \sqrt{\lambda }(S-x)}{\sqrt{\lambda }} + O\left(\frac{e^{|\Im \sqrt{\lambda }|(S-x)}}{\lambda }\right)\right)\left(O\left(\frac{e^{|\Im \sqrt{\lambda }(S+x)|}}{\sqrt{\lambda }}\right)\right)\\
&=& O\left(\frac{e^{2S|\Im \sqrt{\lambda }|}}{\lambda }\right).
\end{eqnarray*}
Thus by Lemma \ref{deltam120} the result holds.  Moreover for $x<0$
\begin{eqnarray*}
&& w_2(\tilde{v}' -v')-v(\tilde{w_2}'-w_2')\\
&=& \left( \frac{\sin \sqrt{\lambda }(x+S)}{\sqrt{\lambda }} + O\left(\frac{e^{|\Im \sqrt{\lambda }|(x+S)}}{\lambda }\right) \right)
\left( O\left(\frac{e^{|\Im \sqrt{\lambda }(S-x)|}}{\sqrt{\lambda }}\right)\right)\\
&-& \left( -m_{22} \frac{\sin \sqrt{\lambda }S}{\sqrt{\lambda }}\cos \sqrt{\lambda }x + m_{11}\frac{\sin \sqrt{\lambda }x}{\sqrt{\lambda }}\cos \sqrt{\lambda }S + O\left(\frac{e^{|\Im \sqrt{\lambda }(S-x)|}}{\lambda }\right) \right) \left(O\left(\frac{e^{|\Im \sqrt{\lambda }|(x+S)}}{\sqrt{\lambda }}\right) \right) \\
&=& O\left(\frac{e^{2S|\Im \sqrt{\lambda }|}}{\lambda }\right).
\end{eqnarray*}
By Lemma \ref{deltam120} the case of $m_{12}=0$ is complete.

Now assume that $m_{12}\neq 0$ and $x>0.$ Then 
\begin{eqnarray*}
&& w_2(\tilde{v}' -v')-v(\tilde{w_2}'-w_2')\\
&=& \left( m_{12} \cos \sqrt{\lambda }S\cos \sqrt{\lambda }x + O\left(\frac{e^{|\Im \sqrt{\lambda }(S+x)|}}{\sqrt{\lambda }}\right) \right)
\left( O\left(\frac{e^{|\Im \sqrt{\lambda }|(S-x)}}{\sqrt{\lambda }}\right) \right)\\
&-& \left(\frac{-\sin \sqrt{\lambda }(S-x)}{\sqrt{\lambda }} + O\left(\frac{e^{|\Im \sqrt{\lambda }|(S-x)}}{\lambda }\right) \right) \left( O\left(e^{|\Im \sqrt{\lambda }(S+x)|}\right)\right) \\
&=& O\left(\frac{e^{2S|\Im \sqrt{\lambda }|}}{\sqrt{\lambda} }\right),
\end{eqnarray*}
from which the result follows using Lemma \ref{deltam12n0}.

Lastly let $m_{12}\neq 0$ and $x<0.$ Then
\begin{eqnarray*}
&& w_2(\tilde{v}' -v')-v(\tilde{w_2}'-w_2')\\
&=& \left ( \frac{\sin \sqrt{\lambda }(x+S)}{\sqrt{\lambda }} + O\left(\frac{e^{|\Im \sqrt{\lambda }|(x+S)}}{\lambda }\right) \right) 
\left( O\left(e^{|\Im \sqrt{\lambda }(S-x)|}\right) \right)\\
&-& \left( -m_{12} \cos \sqrt{\lambda }S \cos \sqrt{\lambda }x + O\left(\frac{e^{|\Im \sqrt{\lambda }(S-x)|}}{\sqrt{\lambda }}\right) \right) \left( O\left(\frac{e^{|\Im \sqrt{\lambda }|(x+S)}}{\sqrt{\lambda }}\right)\right) \\
&=& O\left(\frac{e^{2S|\Im \sqrt{\lambda }|}}{\sqrt{\lambda} }\right).
\end{eqnarray*}
Again using Lemma \ref{deltam12n0} concludes the proof. \qed

%Squaring the contour $\Gamma_k$, we denote this closed contour by $\Gamma_k^2$ 
%and the region enclosed by this contour by $D_k$.  Then $D_1\subset D_2 \subset \dots $
%and $\cup_{k\in\N} D_k=\C$.  So by the maximum modulus principle. 

%%%%%%% bibliography %%%%%%%%%%%%%%%%%%%%%%%%%


\begin{thebibliography}{22}
\bibitem{BBW} {P.A. Binding, P.J. Browne, B.A. Watson,}
      {{Equivalence of inverse Sturm-Liouville problems with boundary conditions rationally dependent on the eigenparameter},
      {\em J. Math. Anal. Appl.}, {\bf 291} (2004), 246--261.}
\bibitem{chadan} {K. Chadan, P.C. Sabatier,}
      {{\em Inverse Problems in Quantum Scattering Theory},
      {Springer-Verlag}, (1977).}
\bibitem{cenw1} {S. Currie, M. Nowaczyk, B.A. Watson,}
		{{Forward scattering on the line with a transfer condition}, {\em Boundary value problems},
      {\bf 2013} (2013) no. 255, 1--14.}      
\bibitem{Freiling} {G. Freiling, V. Yurko},
      {{\em Inverse Sturm-Liouville Problems and their Applications},
      {Nova Science}, (2001).}
\bibitem{Goldberg} {S. Goldberg,}
      {{\em Unbounded Linear Operators, Theory and Applications},
      {McGraw-Hill}, (1966).}      
\bibitem{Gordon} {N.A. Gordon, D.B. Pearson,}
      {{Point Transfer Matrices for the Schr\"odinger Equation: The Algebraic Theory},
      {\em Proceedings of the Royal Society of Edinburgh}, {\bf 129A} (1999), 717--732.}
\bibitem{hald} {O.H. Hald,}
      {{Discontinuous inverse eigenvalue problems},
      {\em Commun. Pure Appl. Math}, {\bf 37} (1984), 539--577.}      
\bibitem{hoch} {H. Hochstadt, B. Lieberman,}
     {{An inverse Sturm-Liouville problem with mixed given data},
      {\em SIAM J. Appl. Math.}, {\bf 34} (1978), 676--680.}      
\bibitem{rost} {R.O. Hryniv,}
      {{Analyticity and uniform stability of the inverse singular Sturm-Liouville spectral problem}, {\em Inverse Problems}, {\bf 27} (2011), 065011 (25pp).}      
\bibitem{Hsieh}{P. Hsieh, Y. Sibuya},
      {{\em Basic Theory of Ordinary Differential Equations},
      {Springer-Verlag}, (1999).}
\bibitem{kes1} {B. Keskin, A.S. Ozkan, N. Yalcin,}
      {{Inverse spectral problems for discontinuous Sturm-Liouville operator with eigenparameter dependent boundary conditions},
      {\em Commun. Fac. Sci. Univ. Ank. Series A1}, {\bf 60} no.1, (2011), 15--25.}     
\bibitem{kes2} {B. Keskin, A. S. Ozkan,}
      {{Spectral problems for Sturm-Liouville operator with boundary and jump conditions linearly dependent on the eigenparameter},
      {\em Inverse Problems in Science and Engineering}, {\bf 20} no.6, (2012), 799--808.}
\bibitem{kes3} {B. Keskin, A.S Ozkan,}
      {{Uniqueness theorems for an impulsive Sturm-Liouville boundary value problem },
      {\em Appl. Math. J. Chinese Univ.}, {\bf 27(4)}, (2012), 428--434.}                               
\bibitem{koba} {M. Kobayashi,}
      {{A uniqueness proof for discontinuous inverse Sturm-Liouville problems with symmetric potentials},
      {\em Inverse Problems}, {\bf 5} (1989), 767--781.}      
\bibitem{levinson} {N. Levinson,}
      {{The inverse Sturm-Liouville problem},
      {\em Matematisk Tidsskrifts}, {\bf B25} (1949a), 25--30.}            
\bibitem{March}{V.A. Mar\v cenko,}
       {{\em Sturm-Liouville Operators and Applications: Revised Edition}
       {AMS}, (2011).}
\bibitem{ramm} {A.G. Ramm,}
      {{One-dimensional inverse scattering and spectral problems},
      {\em CUBO a Math. Journal}, {\bf 6}, N1, (2004), 313--426.}           
\bibitem{Weidmann} {J. Weidmann, }
      {{\em Linear Operators in Hilbert Spaces},
          {Springer-Verlag} (1980).}
\bibitem{Weidmann2} {J. Weidmann, }
      {{\em Spectral Theory of Ordinary Differential Operators},
          {Springer-Verlag} (1980).}
\bibitem{wang} {A. Wang, J. Sun, A. Zettl,}
      {{Two-interval Sturm-Liouville operators in modified Hilbert spaces},
      {\em J. Math. Anal. Appl.}, {\bf 328} (2007), 390--399.}          
\bibitem{willis} {C. Willis,}
      {{Inverse Sturm-Liouville problems with two discontinuities},
      {\em Inverse Problems}, {\bf 1} (1985), 263--289.}
\end{thebibliography}
\end{document}